\documentstyle[12pt,amssymb]{article}
\author{Mohammad~Obiedat} \title{On $J$-orders of
elements of $KO({\Bbb C} P^m)$ }
\newtheorem{thm}{Theorem}[section]
\newtheorem{cor}[thm]{Corollary}

\newtheorem{lemma}[thm]{Lemma} 
\begin{document}
\maketitle \begin{abstract} Let $KO({\Bbb C} P^m)$
be the $KO$-ring of the complex projective space
${\Bbb C} P^m.$ By means of methods of rational
$D$-series \cite{dib1}, a formula for the
$J$-orders of elements of $KO({\Bbb C} P^m)$ is
given.  Explicit formulas are given for computing
the $J$-orders of the canonical generators of
$KO({\Bbb C} P^m)$ and the $J$-order of any
complex line bundle over ${\Bbb C} P^m.$
\end{abstract} \section{Introduction} Let $X$ be a
connected finite-dimensional $CW$ complex.  Let
$KO(X)$ be the $KO$-ring of $X$ and
$\widetilde{K}O(X)$ (resp.  $\widetilde{K}SO(X)$)
be the subgroup of $KO(X)$ of elements (resp.
orientable elements) of virtual dimension zero.
For a real vector bundle $E$ over $X,$ let $S(E)$
be the sphere bundle associated to $E$ with
respect to some inner product on $E$.  Let
$JO(X)=KO(X)/TO(X)$ be the $J$-group of $X,$ where
$TO(X)=\{ E-F \in \widetilde{K}O(X):S(E\oplus n)$
is fibre homotopy equivalent to $S(F\oplus n)$ for
some $ n \in {\Bbb N}\}.$ Then by Adams~\cite{ada}
and Quillen~\cite{qui}, it is shown that
\begin{eqnarray} \label{eq:tox} TO(X)= \{ & x \in
\widetilde{K}SO(X):  \mbox{ there exists } u \in
\widetilde{K}SO(X) \mbox{ such that } \nonumber \\
& \theta_{k}(x)= \frac{\psi^{k}(1+u)}{1+u} \mbox{
in } 1+\widetilde{K}SO(X) \otimes {\Bbb Q}_{k}
\mbox{ for all } k \in {\Bbb N}\} \end{eqnarray}
where $\theta_{k}:\widetilde{K}SO(X)\rightarrow
1+\widetilde{K}SO(X) \otimes {\Bbb Q}_k$ is the
Bott exponential map, $\psi^k$ is the Adams
operation, ${\Bbb Q}_{k}=\{ n/k^m:  n,m \in {\Bbb
Z} \},$ and $ 1+\widetilde{K}SO(X)\otimes {\Bbb
Q}_k$ is the multiplicative group of elements
$1+w$ with $w\in \widetilde{K}SO(X)\otimes {\Bbb
Q}_k.$

 Now, $X$ is connected implies that
$KO(X)=\widetilde{K}O(X)\oplus {\Bbb Z}.$ So,
$JO(X)=\widetilde{J}O(X)\oplus {\Bbb Z}\mbox{
where }\widetilde{J}O(X)=
\widetilde{K}O(X)/TO(X).$ For $x\in KO(X),$ the
$J$-order of $x$ is the order of $x+TO(X)$ in
$JO(X).$ By Atiyah \cite{ati}, $\widetilde{J}O(X)$
is a finite group.  Hence, $x\in KO(X)$ has a
finite $J$-order if and only if
$x\in\widetilde{K}O(X).$ Let $y_m=r \xi_m({\Bbb
C}) -2 $ where $\xi_m({\Bbb C})$ is the complex
Hopf line bundle over the complex projective space
${\Bbb C} P^m.$ Every element of
$\widetilde{K}O({\Bbb C}P^m)$ has the form
$P_{m}(y_m;m_1,\ldots,m_s)=m_1 y_m +\cdots +m_s
y_{m}^s$ for some $m_i\in {\Bbb Z}$ and $s\in
{\Bbb N}.$ The purpose of this paper is to compute
the $J$-order of $P_{m}(y_m;m_1,\ldots,m_s)$ which
we denote by $ b_m(P_m(y_m;m_1,\ldots,m_s)).$ In
addition to their self importance, these orders
are of great importance in geometric topology, for
instance, it is well-known that the Stiefel
fibration $U(n)/U(n-m-1)\longrightarrow S^{2n-1}$
has a cross-section if and only if $n$ is a
multiple of $b_m(P_m(y_m;1,0,\ldots,0)),$ computed
by Adams-Walker~\cite{ada-wal}.
Diba\u{g}~\cite{dib2} has used (\ref{eq:tox}) to
give another proof of Lam's results \cite{lam}, we
rediscover his proof as a special case of Example
$2$ below.

In section~2, we first obtain a useful formula for
 $\theta_{p}(y_{m}^n) $ for $ n=1,\ldots,s.$ Then
 we use (\ref{eq:tox}) and some facts of rational
 $D$-series \cite{dib1} to give, in Theorem~{2.5},
 a formula for $ b_m(P_{m}(y_m;m_1,\ldots,m_s))$ .
 The formula given in Theorem~{2.5} involes many
 difficulties and one can obtain a little
 information about the range of the $J$-orders of
 elements of $\widetilde{K}O({\Bbb C}P^m).$ So
 instead, we use a well-known computations of the
 $J$-order of $y_m$ to obtain, in Theorem~{2.6},
 upper and lower bounds for $
 b_m(P_m(y_m;m_1,\ldots,m_s)).$

There are two important examples in mind, namely
the $J$-orders of the canonical generators of
$\widetilde{K}O ({\Bbb C}P^m)$ and the $J$-orders
of complex line bundles over ${\Bbb C}P^m.$ In
section 3 we first obtain an explicit formula for
the $J$-orders of $y_{m}^{n}$ for $n=2,3,\mbox{ or
}4,$ we conjecture that the formula is liable to
be true for $k=5,\ldots,t.$ Then, we give an
explicit formula for the $J$-order of any complex
line bundle over ${\Bbb C}P^m.$
\section{The $J$-order of
$P_m(y_m;m_1,\ldots,m_s)\in\widetilde{K}O({\Bbb
C}P^m)$} Let $L$ be a non-trivial complex line
bundle over ${\Bbb C}P^m$ .  Then $L \cong
\xi_m({\Bbb C})^n\,$ or $L \cong
\overline{\xi_m({\Bbb C})}^n $ for some $n \in
{\Bbb N},$ where $ \overline{\xi_m({\Bbb C})}$
denotes the conjugate bundle to $\xi_m({\Bbb C})$.
Let $KU({\Bbb C}P^m)$ be the $KU$-ring of ${\Bbb
C}P^m.$ In Lemma~{2.1}, we find the image of
$\xi_m({\Bbb C})^n$ under the realification
homomorphism $$r:KU({\Bbb C}P^m)\rightarrow
KO({\Bbb C}P^m).$$ Note that
$r(\overline{\xi_m({\Bbb C})} ^n)=r(\xi_m({\Bbb
C})^n).$

Recall that from \cite{ada-wal}, $KO({\Bbb C}P^m)$
is a truncated polynomial ring over the integers
generated by $y_m$ with the following relations:
$$ y_{m}^{t+1}=0\, \mbox{ if }\, m=2t; $$ $$
2y_{m}^{2t+1}=0, \, \, y_{m}^{2t+2}=0\, \mbox{ if
}\, m=4t+1;$$ $$ y_{m}^{2t+2}=0 \, \mbox{ if } \,
m=4t+3.  $$ For each $m\in {\Bbb N},$ let $$
d_m=\left\{ \begin{array}{ll} t & \,\mbox{ if }
m=2t \\ 2t+1 & \,\mbox{if } m=4t+1 \mbox{ or }
m=4t+3.  \end{array} \right.  $$

For each $r,s \in {\Bbb N}$ with $r\geq s$, let
$b_{r,s}$ be the coefficient of $(\xi_m({\Bbb
C})^s + \overline{\xi_m({\Bbb C})}^s)$ in
$(\xi_m({\Bbb C}) + \overline{\xi_m({\Bbb C})}
-2)^r.$ Using the fact that $\xi_m({\Bbb C})$
$\overline{\xi_m({\Bbb C})}=1,$ we easily see
that:  \begin{equation} \label{eq:brs}
b_{r,s}=(-2)^{r-s}{r\choose s}+ \sum_{{j=s+1}
\atop {j=s+2s_j}}^{r} (-2)^{r-j} {r \choose j} {j
\choose s_j}.  \end{equation} Further, $b_{r,0} $
is the constant term of $(\xi_m({\Bbb C}) +
\overline{\xi_m({\Bbb C})} -2)^r.$\\

Let $n \in {\Bbb N}.$ For each $s=1,\ldots,n,$
define $d_{n,s}$ by the recurrence relation
$d_{n,n} =1$ and for $s=n-1,n-2,\ldots,1$
\begin{equation} \label{eq:dns}
d_{n,s}=-(d_{n,s+1} b_{s+1,s}+d_{n,s+2}
b_{s+2,s}+\cdots +d_{n,n} b_{n,s}).
\end{equation} \textbf{Convention.}  Let $f\in
{\Bbb Z} [\![ y_m]\!],$ where $ {\Bbb Z} [\![
y_m]\!]$ is the ring of formal power series with
coefficients in ${\Bbb Z}.$ If we consider $f$ as
an element of $KO({\Bbb C}P^m),$ then we
implicitly mean that $f\,(\mbox{mod }
y_{m}^{d_{m}+1}).$
\begin{lemma} Let $n,m \in {\Bbb N}.$ Then
\begin{enumerate} \item[(i)] $r(\xi_m({\Bbb C})
^n)= d_{n,1} y_{m} + d_{n,2}
y_{m}^{2}+\cdots+d_{n,n} y_{m}^{n}+2.$ \item[(ii)]
$\psi^n(y_m)=d_{n,1} y_{m} + d_{n,2}
y_{m}^{2}+\cdots+d_{n,n} y_{m}^{n},$ and for
$k=2,\ldots,d_m,$ $$\psi^n(y_{m}^{k})=d_{kn,1}
y_{m} + \cdots+d_{kn,kn} y_{m}^{kn}-d_{k,1}
\psi^n(y_{m})-\cdots-d_{k,k-1}\psi^n(y_{m}^{k-1}).$$
\item[(iii)] $d_{n,1}=n^2 .$ \end{enumerate}
\end{lemma} \textbf{Proof.}  (i) Let $c:KO({\Bbb
C}P^m)\rightarrow KU({\Bbb C}P^m)$ be the
complexification homomorphism.  $cr(\xi_m({\Bbb
C})^n)=\xi_m({\Bbb C})^n+\overline{\xi_m({\Bbb
C})}^n.$ On the other hand, by
(\ref{eq:brs}),(\ref{eq:dns}) and the fact that
$y_{m}^n = 0$ for $n >d_m,$ we have $$c(d_{n,1}
y_{m} +\cdots+ d_{n,n} y_{m}^{n}+2)=\xi_m({\Bbb
C})^n+\overline{\xi_m({\Bbb C})}^n.$$ Using the
fact that $c$ is a monomorphism for $m=2t$ and
$m=4t+3,$ we get $$r(\xi_m({\Bbb C}) ^n)= d_{n,1}
y_{m} +\cdots+ d_{n,n} y_{m}^{n}+2.$$ To prove the
case $m=4t+1,$ let $i:{\Bbb C} P^{4t+1}
\rightarrow {\Bbb C} P^{4t+2}$ be the inclusion
map.  Then $i^*:KO({\Bbb C} P^{4t+2}) \rightarrow
KO({\Bbb C} P^{4t+1})$ is an epimorphism and maps
$r(\xi_{4t+2}({\Bbb C})^n)$ to $r(\xi_{4t+1}({\Bbb
C})^n).$ Hence, $$r(\xi_{4t+1}({\Bbb C})^n)= i^*(
r(\xi_{4t+2}({\Bbb
C})^n))=i^*(d_{n,1}y_{4t+2}+\cdots+d_{n,n}y_{4t+2}^n+2)$$
$$=d_{n,1}y_{4t+1}+d_{n,2}y_{4t+1}^2+\cdots+d_{n,n}y_{4t+1}^n+2.
$$ (ii) Let $l\in \{ 1,\ldots,d_m\}.$ By (i),
$$d_{ln,1} y_{m} + \cdots+d_{ln,ln}
y_{m}^{ln}=r(\xi_m({\Bbb C})^{ln})-2
=r\psi^n(\xi_m({\Bbb C})^{l})-2$$
$$=\psi^n(r\xi_m({\Bbb C})^{l})-2= \psi^n(d_{l,1}
y_{m} + \cdots+d_{l,l} y_{m}^{l}).$$ The result
follows.\\ (iv) $d_{n,1}$ is the constant term of
$$\frac{cr(\xi_m({\Bbb C})^n)-2} {\xi_m({\Bbb
C})+\overline{\xi_m({\Bbb C})}-2}=
\frac{\xi_m({\Bbb C})^n+\overline{\xi_m({\Bbb
C})}^n-2}{\xi_m({\Bbb C})+ \overline{\xi_m({\Bbb
C})}-2}=$$ $$(\xi_m({\Bbb C})^{n-1}+\xi_m({\Bbb
C})^{n-2}+\cdots+\xi_m({\Bbb C})+1)(
(\overline{\xi_m({\Bbb
C})}^{n-1}+\overline{\xi_m({\Bbb C})}^{n-2}+
\cdots+\overline{\xi_m({\Bbb C})}+1).$$ Hence,
$$d_{n,1}=n+2(n-1)+2(n-2)+ \cdots +2$$
$$=n+2((n-1)+(n-2)+\cdots+1)=n+2(\frac{n(n-1)}{2})=n^2.$$
This completes the proof of Lemma~{2.1.}

Now, we use the above lemma to find
$\theta_{p}(y_{m}^n) \in 1+\widetilde{K}O({\Bbb
C}P^m) \otimes{\Bbb Q}_p.$ For each $n,m,p\in
{\Bbb N}$ with $n\leq d_m,$ let $$A(p;n,m)=\left(
\frac{\psi^p(d_{n,1}+\cdots+
d_{n,n}y_{m}^{n-1})}{d_{n,1}+\cdots+d_{n,n}y_{m}^{n-1}}
\right)^{\frac{1}{2}}\in 1+\widetilde{K}O({\Bbb
C}P^m)\otimes{\Bbb Q}_p,$$ and for $n\geq 2,$ let
$$B(p;n,m)\!=\left(\theta_p(y_{m})^{d_{n,1}-1}
\theta_p(y_{m}^2)^{d_{n,2}} \ldots
\theta_p(y_{m}^{n-1})^{d_{n,n-1}}\right)^{-1}\!\in
\!1+ \widetilde{K}O({\Bbb C}P^m) \otimes{\Bbb
Q}_p.$$
\begin{thm}
Let $p\geq 2$ and $m=2t$ for some $t \geq 1$.
Then \begin{enumerate} \item[(i)]
$\theta_p(y_m)=\left(\frac{\psi^p(y_m)}{p^2y_m}\right)^\frac{1}{2}.$
\item[(ii)] $\theta_p(y_{m}^n)=A(p;n,m)B(p;n,m),$
for each $2\leq n \leq t.$ \end{enumerate}
\end{thm} \textbf{Proof.}  (i) This is Lemma 5.4
of \cite{mps}.\\ (ii) If $\eta$ is a complex
$4n$-dimensional vector bundle over a finite $CW$
complex $X$ such that $\wedge^{4n}\eta=1,$ then
$c\theta_p(r\eta)=\theta_p(\eta).$ Let $\eta
=2\xi_m({\Bbb C})^n+2 \overline{\xi_m({\Bbb
C})}^n.$ Then $\eta$ is a 4-dimensional complex
vector bundle over ${\Bbb C}P^m$ with
$$\wedge^4(\eta)=\wedge^2 (2\xi_m({\Bbb C})^n)
\wedge^2(2\overline{\xi_m({\Bbb C})}^n)=1.$$
Hence, $ c \theta_p(r\eta)=\theta_p(\eta).$ By
Lemma~ 2.1, $$\xi_m({\Bbb
C})^n+\overline{\xi_m({\Bbb C})}^n-2=
cr(\xi_m({\Bbb
C})^n)-2=c(d_{n,1}y_m+\cdots+d_{n,n} y_{m}^n).$$
Also, $$r\eta =2r(cr(\xi_m({\Bbb
C})^n))=2rc(r(\xi_m({\Bbb C})^n))$$
$$=4r(\xi_m({\Bbb C})^n)=
4d_{n,1}y_m+\cdots+4d_{n,n} y_{m}^n+8.$$ Thus,
$$c\theta_p(r\eta)=c((p \theta_p(y_{m})^{d_{n,1}}
\theta_p(y_{m}^2)^{d_{n,2}} \ldots
\theta_p(y_{m}^{n-1})^{d_{n,n-1}}
\theta_p(y_{m}^n))^4).$$ On the other hand,
$$\theta_p(\eta)= \left( \frac{\xi_m({\Bbb
C})^{np}+ \overline{\xi_m({\Bbb
C})}^{np}-2}{\xi_m({\Bbb C})^n+
\overline{\xi_m({\Bbb C})}^n-2 }
\right)^2=\left(\frac{\psi^p(\xi_m({\Bbb C})^n+
\overline{\xi_m({\Bbb C})}^n-2)}{\xi_m({\Bbb
C})^n+ \overline{\xi_m({\Bbb C})}^n-2 } \right)^2
$$ $$= \left( \frac{c(\psi^p(y_m))
c(\psi^p(d_{n,1}+d_{n,2}y_m+\cdots+
d_{n,n}y_{m}^{n-1}))}{c(y_m)
c(d_{n,1}+d_{n,2}y_m+\cdots+ d_{n,n}y_{m}^{n-1})}
\right)^2.$$ By (i),
$$c\left(\frac{\psi^p(y_m)}{y_m}\right)=c(p^2\theta_p(y_m)^2).$$
Hence, $$\theta_p(\eta)=c \left( (p
\theta_p(y_m))^4(\frac{\psi^p(d_{n,1}+d_{n,2}y_m+\cdots+
d_{n,n}y_{m}^{n-1})}{d_{n,1}+d_{n,2}y_m+\cdots+
d_{n,n}y_{m}^{n-1}})^2 \right).$$ Now, c is a
monomorphism implies that $$p
\theta_p(y_{m})^{d_{n,1}}
\theta_p(y_{m}^2)^{d_{n,2}} \ldots
\theta_p(y_{m}^{n-1})^{d_{n,n-1}}
\theta_p(y_{m}^n)$$ $$=p\theta_p(y_m) \left(
\frac{\psi^p(d_{n,1}+d_{n,2}y_m+\cdots+d_{n,n}y_{m}^{n-1})}{d_{n,1}+d_{n,2}y_m+\cdots+
d_{n,n}y_{m}^{n-1}} \right)^\frac{1}{2}.$$ Hence,
$$\theta_p(y_{m}^n)B(p;n,m)^{-1}=A(p;n,m).$$ This
completes the proof of Theorem~{2.2.}

REMARK.  By Using a method similar to that used in
proving Lemma~{2.1} when $m=4t+1$, we easily
obtain a similar formula for $\theta_p(y_{m}^n)$
when m is an odd integer.

\begin{cor} $\theta_p \circ \psi^n =\psi^n \circ
\theta_p$ on $\widetilde{K}O({\Bbb C}P^m)$ for all
$m,n,p \in {\Bbb N}.$ \end{cor} \textbf{Proof.}
Clearly, we only need to show that $\theta_p \circ
\psi^n(y_{m}^{k}) =\psi^n \circ
\theta_p(y_{m}^{k})$ for each $k=1,\ldots,d_m.$ By
induction on $k,$ if $k=1$ then by Theorem~{2.2}
and Lemma~{2.1} (ii), $$\theta_p (
\psi^n(y_{m}))=\theta_p(d_{n,1}y_m+\cdots+d_{n,n}y_{m}^{n})$$
$$= \left(\frac{\psi^p(d_{n,1}+d_{n,2}y_{m}\cdots+
d_{n,n}y_{m}^{n-1})}{d_{n,1}+d_{n,2}y_m+\cdots+d_{n,n}y_{m}^{n-1}}
\right)^\frac{1}{2}\theta_p(y_m)$$
$$=\left(\frac{\psi^p(\psi^n(y_{m}))}{p^2\psi^n(y_m)}
\right)^\frac{1}{2}
=\psi^n\left(\frac{\psi^p(y_{m})}{p^2y_m}
\right)^\frac{1}{2}=\psi^n(\theta_p(y_m)).$$ Now,
suppose $\theta_p ( \psi^n(y_{m}^{l}) )=\psi^n (
\theta_p(y_{m}^{l}))$ for $l=1,\ldots,k-1.$ Then
$$\theta_p ( \psi^n(y_{m}^{k}) )=\theta_p(d_{kn,1}
y_{m} + \cdots+d_{kn,kn} y_{m}^{kn}-d_{k,1}
\psi^n(y_{m})-\cdots-d_{k,k-1}\psi^n(y_{m}^{k-1}))$$
$$=\theta_p(y_m)^{d_{kn,1}}\cdots\theta_p(y_{m}^{kn})^{d_{kn,kn}}
\theta_p(\psi^n(y_m))^{-d_{k,1}}\cdots\theta_p(\psi^n(y_{m}^{k-1}))^{-d_{k,k-1}}.$$
On the other hand,
$\psi^n(\theta_p(y_{m}^{k}))=\psi^n(A(p;k,m)B(p;k,m)).$
So, we only need to show that
$$\theta_p(y_m)^{d_{kn,1}}\cdots\theta_p(y_{m}^{kn})^{d_{kn,kn}}=
A(p;k,m)\psi^n(\theta_p(y_m)).$$
$$\theta_p(y_m)^{d_{kn,1}}\cdots\theta_p(y_{m}^{kn})^{d_{kn,kn}}=
\left(\frac{\psi^p(d_{kn,1}+d_{kn,2}y_{m}^{2}
\cdots+d_{kn,kn}y_{m}^{kn-1})}{d_{kn,1}+d_{kn,2}y_m+\cdots+
d_{kn,kn}y_{m}^{kn-1}}
\right)^\frac{1}{2}\theta_p(y_m)$$
$$=\left(\frac{\psi^p(d_{kn,1}y_m+d_{kn,2}y_{m}^{2}\cdots+
d_{kn,kn}y_{m}^{kn})}{p^2(d_{kn,1}y_m+d_{kn,2}y_{m}^{2}+\cdots+
d_{kn,kn}y_{m}^{kn})} \right)^\frac{1}{2}
=\left(\frac{\psi^p(\psi^{kn}(y_{m}))}{p^2\psi^{kn}(y_m)}
\right)^\frac{1}{2}$$
$$=\left(\psi^n(\frac{\psi^p(\psi^k(y_{m}))}{p^2\psi^k(y_m)})
\right)^\frac{1}{2}=
A(p;k,m)\psi^n(\theta_p(y_m)).$$ The result
follows.

Let
$s=d_m$ and
$P_m(y_m;m_1,\ldots,m_s)=m_1y_m+\cdots+m_sy_{m}^{s}\in
\widetilde{K}O({\Bbb C}P^m).$ Let
$b_m(P_m(y_m;m_1,\ldots,m_s))$ be the $J$-order of
$P_m(y_m;m_1,\ldots,m_s),$ that is the order of
$P_m(y_m;m_1,\ldots,m_s)+TO({\Bbb C}P^m)$ in
$\widetilde{J}O({\Bbb C}P^m).$ By using the same
method of Proposition~{5.7} of \cite{mps}, we
easily see that
$$b_{4t+3}(P_{4t+3}(y_{4t+3};m_1,\ldots,m_{2t+1}))=
b_{4t+2}(P_{4t+2}(y_{4t+2};m_1,\\ \ldots,
m_{2t+1})),$$ and
$b_{4t+1}(P_{4t+1}(y_{4t+1};m_1,\ldots,m_{2t+1})$
$$=\left\{ \begin{array}{ll}
b_{4t}(P_{4t}(y_{4t};m_1,\ldots,m_{2t}))& \mbox{if
} m_{2t+1}=0\\
lcm\{b_{4t}(P_{4t}(y_{4t};m_1,\ldots,m_{2t})),2 \}
& \mbox{if } m_{2t+1}=1.  \end{array} \right.  $$
Therefore, in the remainder of this paper we shall
assume that $m=2t$ for some $t\geq 1$, unless
othewise indicated.

To compute $b_m(P_m(y_m;m_1,\ldots,m_t))$ , we use
the notion of rational $D$-series introduced and
devloped in \cite{dib1}.  Let ${\Bbb Q} [\![x]\!]$
be the ring of formal power series with
coefficients in ${\Bbb Q}$ and ${\Bbb Q}^*
[\![x]\!]=\{f(x) \in {\Bbb Q} [\![x]\!]:f(0)=\pm
1\}.$ Let $$f(x)=\pm 1+\sum_{i\geq 1} a_i x^i \in
{\Bbb Q}^* [\![x]\!].$$ For each $k \geq 1,$ let
$e_k(f)$ denote the smallest positive integer
$e_k$ such that $(f(x))^{e_k} \in {\Bbb
Z}[\![x]\!]  (\mbox{mod }x^{k+1}).$ Let $S_k(f)$
be the set of all primes dividing the denominators
of the coefficients $a_i$ for $i=1,\ldots,k.$ For
a rational number $q$, let $\nu_p(q)$ be the
exponent of p in the prime factorization of $q$
and let $D(q)$ be the denominator of $q$ in its
lowest term.  For convenience, we assume that
$\nu_p(0)=-\infty $ and $D(0)=1.$ It follows from
Lemma~{1.3} \cite{dib1} that $p \in S_k(f)$ if and
only if $\nu _p(e_k(f)) > 0.$

For each prime $p$, let $\alpha_p,\beta_p \in
{\Bbb Z}^+,$ and let $\alpha\!=\!\!
(\!\alpha_2,\alpha_3,\alpha_5,\ldots\!)$, and
$\beta\!=\!\!(\!\beta_2,\beta_3,\beta_5,\ldots\!).$
A series $f=\pm 1+\sum_{i\geq 1} a_i x^i\in {\Bbb
Q}^* [\![x]\!]$ is called a rational $D$-series of
type $(\alpha,\beta)$ if
$\nu_p(a_{\alpha_p})=-\beta_p$ and $ \nu_p(a_k)
\geq -\beta_p[\frac{k}{\alpha_p}]$ for each $k >
\alpha_p$ with $a_k\neq 0.$ A rational $D$-series
$f$ is called strict at a prime $p$ if
$\nu_p(a_k)>-\beta_p [\frac{k}{\alpha_p}]$ for
each $k > \alpha_p$ with $a_k\neq 0.$ If $f$ is
strict at $p$ then, by Theorem~{3.5} of
\cite{dib1}, $$\nu_p(e_k(f))=\max\{ \beta_p
r+\nu_p(r):0 \leq r \leq [\frac{k}{\alpha_p}]
\}.$$ With these facts on hand, we return to our
problem.

Let ${\Bbb Z}_{(p)}=\{\frac{r}{s}:r,s\in {\Bbb Z}
\mbox{ with } \nu_p(s)=0\}$ be the localization of
${\Bbb Z}$ at $p.$ The following lemma is (5.2)
and Lemma 5.5 of \cite{mps} with minor changes.
\begin{lemma} \begin{enumerate} \item[(i)] Let
$1+u$ be an element of ${\Bbb
Q}^*[\![y_m]\!](\mbox{mod }y^{t+1}).$ Then
$$\frac{\psi^p(1+u)}{1+u}\in {\Bbb
Z}_{(p)}^*[\![y_m]\!](\mbox{mod }y_{m}^{t+1})$$ if
and only if $u \in {\Bbb
Z}_{(p)}[\![y_m]\!](\mbox{mod }y^{t+1})$ with
$u(0)=0.$ \item[(ii)] $(\theta_p(y_m))^h \in {\Bbb
Z}_{(p)}^*[\![y_m]\!]  (\mbox{mod }y_{m}^{t+1}) $
if and only if $$\nu_p(h) \geq
\nu_p(b_m(y_m))=\max \{ s+\nu_p(s):0 \leq s \leq
[\frac{m}{p-1}] \}.$$ \end{enumerate} \end{lemma}
Now, we
compute $b_m(P_m(y_m;m_1,\ldots,m_t)).$ For each
$n=1,\ldots,t,$ let
$\theta_p(y_{m}^n)=1+\alpha_{n,1} y_m
+\cdots+\alpha_{n,t} y_{m}^t$ where $\alpha_{n,i}$
is the coefficient of $y_{m}^i$ given by Theorem
2.2.  According to (\ref{eq:tox})
$b_m(P_m(y_m;m_1,\ldots,m_t))$ is the smallest
positive integer $h$ such that $$\theta_p(h
P_m(y_m;m_1,\ldots,m_t))=\frac{\psi^p(1+u)}{1+u}\mbox{
in }\,\, \widetilde{K}O({\Bbb C}P^m)\otimes{\Bbb
Q}$$ for some $u\in \widetilde{K}O({\Bbb C}P^m)$
and all primes $p.$ Let
$$\beta_m(y_m;m_1,\ldots,m_t)=\theta_p(P_m(y_m;m_1,\ldots,m_t))=
\theta_p(y_{m})^{m_1} \ldots
\theta_p(y_{m}^t)^{m_t} $$
$$=1+\alpha_1(m_1,\ldots,m_t)y_m+\cdots+
\alpha_t(m_1,\ldots,m_t)y_{m}^t,$$ for some
$\alpha_i(m_1,\ldots,m_t)\in {\Bbb Q}.$ Then
$$\theta_p(h
P_m(y_m;m_1,\ldots,m_t))=\beta_m(y_m;m_1,\ldots,m_t)^h
.$$ Let
$$\beta_m(y_m;m_1,\ldots,m_t)^h=\frac{\psi^p(1+u)}{1+u}$$
for some $u \in\widetilde{K}O({\Bbb C}P^m)$ .
Then $\beta_m(y_m;m_1,\ldots,m_t)^h$ has integer
coefficients.  Hence $$\nu_p(h)\geq
\nu_p(e_t(\beta_m(y_m;m_1,\ldots,m_t))).$$ On the
other hand, if
$b=e_t(\beta_m(y_m;m_1,\ldots,m_t))$ then
$$\beta_m(y_m;m_1,\ldots,m_t)^b=(\frac{\psi^p(1+u)}{1+u})^{\frac{b}{h}}
=\frac{\psi^p(1+w)}{1+w}$$ for some $w\in {\Bbb
Q}^*[y_m](\mbox{mod }y_{m}^{t+1}).$
$\beta_m(y_m;m_1,\ldots,m_t)^b\in {\Bbb
Z}^*[y_m](\mbox{mod }y_{m}^{t+1})$ implies that
$$\frac{\psi^p(1+w)}{1+w}\in {\Bbb
Z}^*[y_m](\mbox{mod }y_{m}^{t+1}).$$ So, by Lemma
2.4 $w\in \widetilde{K}O({\Bbb C}P^m)$ and hence
$\nu_p(h)\leq\nu_p(b).$ By Corollary~{1.3} of
\cite{dib1}, $e_t(\beta_m(y_m;m_1,\ldots,m_t))=
D(b_1) \cdots D(b_t)$ where
$b_1=\alpha_1(m_1,\ldots,m_t)$ and for
$k=1,\dots,t,$ $b_k$ is the coefficient of
$y_{m}^k$ in $\beta_m(y_m;m_1,\ldots,m_t)^{D(b_1)
\ldots D(b_{k-1})}.$ So far, we have proved:
\begin{thm} $b_m(P_m(y_m;m_1,\ldots,m_t))=D(b_1)
\cdots D(b_t).$ \end{thm}
Although Theorem~{2.5} gives
$b_m(P_m(y_m;m_1,\ldots,m_t))$ by a formula, it is
difficult to use this formula to find
$b_m(P_m(y_m;m_1,\ldots,m_t))$ for specific values
of $m_1,\ldots,m_t,$ because one needs first to
find the coefficients of $y_{m}^k$ in
$\theta_p(y_{m})^{m_1} \ldots
\theta_p(y_{m}^t)^{m_t}$ and then to find
$e_t(\beta_m(y_m;m_1,\ldots,m_t))$ which involves
tedious calculations.  So, alternatively, we next
try to obtain information about
$b_m(P_m(y_m;m_1,\ldots,m_t))$ by using what we
know about $b_m(y_m).$

By Theorem 2.2, we directly obtain
 $$\theta_p(y_{m}^k)=\frac{\alpha_k}
 {\theta_p(y_m)^{N_k}}\, \mbox{ where }
 \alpha_2=(\frac{\psi^p(d_{2,1}
 +y_m)}{d_{2,1}+y_m})^{\frac{1}{2}},\,
 N_2=d_{2,1}-1$$ and for $k=3,\ldots,t, $
 \begin{eqnarray} \label{eq:alk} \alpha_k =\left(
 \frac{\psi^p(d_{k,1}+d_{k,2}y_m+\cdots+
 d_{k,k}y_{m}^{k-1})}{d_{k,1}+d_{k,2}y_m+\cdots+
 d_{k,k}y_{m}^{k-1}}
 \right)^{\frac{1}{2}}\alpha_{2}^{-d_{k,2}}
 \alpha_{3}^{-d_{k,3}}\ldots
 \alpha_{k-1}^{-d_{k,k-1}},\nonumber \\
 N_k=(d_{k,1}-1)-N_2
 d_{k,2}-\cdots-N_{k-1}d_{k,k-1}.  \end{eqnarray}
 For each $m_1,\ldots,m_t\in {\Bbb Z},$ let
 $$E(m_1,\ldots,m_t)= \mbox{ lcm}
 \{e_t(\alpha_{2}^{m_2}),\ldots,e_t(\alpha_{t}^{m_t})\},$$
 $$N(m_1,\ldots,m_t)=m_1-m_2N_2-\cdots-m_tN_t,$$
 $$L(p;m_1,\ldots,m_t)=\nu_p(b_m(y_m))-\nu_p(N(m_1,\ldots,m_t))-
 \nu_p(E(m_1,\ldots,m_t)),$$ and
 $$U(p;m_1,\ldots,m_t)=\max\{\nu_p(b_m(y_m))-\nu_p(N(m_1,\ldots,m_t)),
 \nu_p(E(m_1,\ldots,m_t))\}.$$
 \begin{thm} Let
 $P_m(y_m;m_1,\ldots,m_t)=m_1y_m+\cdots+m_ty_{m}^{t}\in
 \widetilde{K}O({\Bbb C}P^m).$ Then
 $$L(p;m_1,\ldots,m_t) \leq \nu_p(b_m(
 P_m(y_m;m_1,\ldots,m_t)) \leq
 U(p;m_1,\ldots,m_t).$$ \end{thm} \textbf{Proof.}
 Let $ h=b_m(P_m(y_m;m_1,\ldots,m_t)).$ Then
 $$\theta_p(m_1y_m+\cdots+m_ty_{m}^t)^h=\frac{\psi^p(1+u)}{1+u}$$
 for some $u \in\widetilde{K}O({\Bbb C}P^m)$ and
 all primes $p.$ So,
 $$\theta_p(y_m)^{N(m_1,\ldots,m_t)h}=
 \frac{\psi^p(1+u)}{1+u}\alpha_{2}^{-m_2h}
 \ldots \alpha_{t}^{-m_th}.$$ Thus,
 $$\theta_p(y_m)^{N(m_1,\ldots,m_t)E(m_1,\ldots,m_t)h}$$
 has integer coefficients.  Hence, by Lemma~{2.4}
 (ii) $$\nu_p(N(m_1,\ldots,m_t))
 +\nu_p(E(m_1,\ldots,m_t))+\nu_p(h) \geq
 \nu_p(b_m(y_m)),$$ namely $\nu_p(h) \geq
 L(p;m_1,\ldots,m_t).$ On the other hand, let $b
 \in {\Bbb N}$ such that
 $\nu_p(b)=U(p;m_1,\ldots,m_t).$ Then
 $$\theta_p(P_m(y_m;m_1,\ldots,m_t))^b=
 (\frac{\psi^p(1+u)}{1+u})^{\frac{b}{h}}.$$
 So, $$\theta_p(y_m)^{N(m_1,\ldots,m_t)b}=
 (\frac{\psi^p(1+u)}{1+u} )^\frac{b}{h}
 \alpha_{2}^{-m_2b} \ldots \alpha_{t}^{-m_tb}.$$
 Let
 $$\left(\frac{\psi^p(1+u)}{1+u}\right)^{\frac{b}{h}}=
 \frac{\psi^p(1+w)}{1+w}$$
 for some $w\in {\Bbb Q}^*[\![y_m]\!](\mbox{mod
 }y_{m}^{t+1})$ with $w(0)=0.$ Now,
 $$\alpha_{2}^{-m_2b} \ldots
 \alpha_{t}^{-m_tb}\mbox{ and }
 \theta_p(y_m)^{N(m_1,\ldots,m_t)b} $$ have
 integer coefficients.  Hence,
 $$\frac{\psi^p(1+w)}{1+w} $$ has integer
 coefficients, which implies that $w\in {\Bbb
 Z}^*[\![y_m]\!](\mbox{mod }y_{m}^{t+1}).$ Hence,
 $$\nu_p(h) \leq \nu_p(b).$$ This completes the
 proof of Theorem~{2.6.}

\begin{cor} Let
$P_m(y_m;m_1,\ldots,m_t)=m_1y_{m}+\cdots+m_ty_{m}^t
\in \widetilde{K}O({\Bbb C}P^m).$ Let $s$ be the
smallest positive integer such that $m_i=0$ for
all $i>s,$ if $m_t \neq 0,$ let $s=t$.  Then
$$\nu_p(b_m(P_m(y_m;m_1,\ldots,m_t))=
\max\{\nu_p(b_m(y_m))-\nu_p(N(m_1,\ldots,m_t)),
0\}$$ for all $p>s.$ \end{cor} \textbf{Proof.}
Using Lemma {2.1} (iii), we easily see that
$S_t(\alpha_{k}^{m_k}) \subseteq \{
2,3,\ldots,k\}$ for each $k=2,\ldots,t.$ So, if
$p>s$ then $\nu_p(E(m_1,\ldots,m_t))=0.$ Now, the
result follows from Theorem~{2.6}.

\section{Two Important Examples} Let $m=2t.$
Then$$\widetilde{J}O({\Bbb C}P^m)=< \alpha_1=y_m+
TO({\Bbb C}P^m),\ldots,\alpha_t=y_{m}^t+TO({\Bbb
C}P^m)>.$$ In Example~1, we give a simple formula
for the $J$-orders of $\alpha_2,\alpha_3,$ and
$\alpha_4.$ We conjecture that our formula is
liable to be true for the $J$-orders of
$\alpha_k,\, k=5,\ldots,t.$

Let $L \cong \xi_m({\Bbb C})^n$ for some $n \in
{\Bbb N}.$ By the $J$-order of $L$ we mean the
order of $r\xi_m({\Bbb C})^n-2+TO({\Bbb C}P^m)$ in
$JO({\Bbb C}P^m).$ Lam \cite{lam} has used complex
K-theory to find the $J$-order of $L$ when $n$ is
a prime power.  Also, Diba\u{g} \cite{dib2} has
used (\ref{eq:tox}) to give another proof of Lam's
results.  In Example~2, we give a simple formula
for the $J$-order of $L$ for each $n \in {\Bbb
N}.$\\
EXAMPLE~{1.}  For each $k=2,\ldots,t,$ the
$J$-order of $\alpha_k=y_{m}^k+TO({\Bbb C}P^m)$ is
$$b_m(P_m(y_m;0,\ldots,0,m_k=1,0,\ldots,0)).$$ So,
by Corollary~{2.7}, if $p>k$ then
$$\nu_p(b_m(y_{m}^k))=\nu_p(b_m(y_m))-\nu_p(N_k).$$
The first few values of $N_k$ are
$N_2=3,\,N_3=-2.5,\,N_4=5.7,\,N_5=-2.3^2.7,\,N_6=2.3.7.11.$
By inspection, if $p>k$ then
$\nu_p(N_k)=[\frac{2(k-1)}{p-1}]$ for each
$k=2,\ldots,6.$

In \cite{obi}, we proved that if $p=2,$ or $3,$
then $$\nu_p(b_m(y_{m}^k))=\max
\{r-[\frac{2(k-1)}{p-1}]+\nu_p(r):
[\frac{2k}{p-1}]\leq r\leq [\frac{m}{p-1}] \}$$
for each $k=1,\ldots,t.$ So, we have:
\begin{thm}
If $k=2,3, \mbox{ or } 4,$ then
$$\nu_p(b_m(y_{m}^k))=\max
\{r-[\frac{2(k-1)}{p-1}]+\nu_p(r):
[\frac{2k}{p-1}]\leq r\leq [\frac{m}{p-1}] \}$$
for each $p\geq 2.$ Further, this formula is true
if $p=2,3$ and $k=5,\dots,t$ or if $k=5, \mbox{ or} 6 
\mbox{ and } p>5.$ \end{thm}
REMARK.  From the above discussion, it is
reasonable to conjecture that :  If $p$ is any
prime number and $k=5,\ldots,t,$ then
$$\nu_p(b_m(y_{m}^k))=\max
\{r-[\frac{2(k-1)}{p-1}]+\nu_p(r):
[\frac{2k}{p-1}]\leq r\leq [\frac{m}{p-1}] \}.$$
Now, we
compute the $J$-order of any complex line bundle
over ${\Bbb C}P^m.$\\ EXAMPLE~{2.}  Let $n\in
{\Bbb N}.$ Then the $J$-order of $\xi_m({\Bbb
C})^n$ is $b_m(r\xi_m({\Bbb C})^n-2).$ By
Theorem~{2.5}, $\nu_p(b_m(r\xi_m({\Bbb C})-2))=
\nu_p(e_t(\theta_p(r\xi_m({\Bbb C})^n-2))).$
$$r\xi_m({\Bbb C})^n-2=r\psi^n(\xi_m({\Bbb C})-1)=
\psi^n(r\xi_m({\Bbb C})-2)=\psi^n(y_m).$$ So, by
Corollary~{2.3}, $$\theta_p(r\xi_m({\Bbb
C})^n-2)=\theta_p(\psi^n(y_m))=
\psi^n(\theta_p(y_m)).$$ Let $n=p_{1}^{r_1} \ldots
p_{s}^{r_s} p^d $ where $r_i >0$ for
$i=1,\ldots,s,$ and $d\geq 0.$ Let $p=2q+1$ be any
odd prime number.  By, Lemma~{5.4} \cite{mps},
$$\theta_p(y_m)=1+\sum_{j=1}^{q-1}m_jy_{m}^{j}+\frac{1}{p}y_{m}^{q}$$
where $m_j\in {\Bbb Z}$ with $\nu_p(m_j)=0$ for
$j=1,\ldots,q-1.$ Hence,
$$\psi^n(\theta_p(y_m))=1+\sum_{j=1}^{q-1}m_j\psi^n(y_{m})^{j}+
\frac{1}{p}\psi^n(y_{m})^{q}.$$ Using Lemma~{3.6}
of \cite{mps}, we easily obtain that
$$\psi^{p^d}(y_m)=\sum_{j=1}^{p^d-1}n_jy_{m}^{j}+y_{m}^{p^d}$$
with $\nu_p(n_j)>0$ for $j=1,\ldots,p^d-1.$ Now,
using the fact that Adams operations are ring
homomorphisms with
$\psi^{l_1}\circ\psi^{l_2}=\psi^{l_1l_2}$ for each
$l_1,l_2\in{\Bbb Z},$ we get
$\psi^n(y_m)=\sum_{j=1}^{n}a_jy_{m}^{j}$ where
$a_j\in{\Bbb Z}$ such that $\nu_p(a_j)=0$ if $j$
is a multiple of $p^d,$ and $\nu_p(a_j)\geq 1$ for
all other values of $j.$ Hence,
$$\psi^n(\theta_p(y_m))=1+\sum_{j=1}^{nq}b_jy_{m}^{j}$$
such that $\nu_p(b_j)=-1$ if $j$ is a multiple of
$p^dq,$ $\nu_p(b_j)\geq 1$ for $j<p^dq,$ and
$\nu_p(b_j)\geq 0$ for all other values of $j.$
Hence, $\theta_p(\psi^n(y_m))$ is a strict
D-series at $p$ of type $(\alpha,\beta)$ where $$
\alpha_{p^{\prime}}=\left\{ \begin{array}{ll} p^dq
& p^{\prime}=p \\ \infty & p^{\prime} \neq p
\end{array} \right.  \hspace{2cm}, \,
\beta_{p^{\prime}}=\left\{ \begin{array} {ll} 1 &
p^{\prime}=p \\ \infty & p^{\prime} \neq p.
\end{array}\right.  $$ Hence,
$$\nu_p(b_m(r\xi_m({\Bbb
C})^n-2))=\max\{r+\nu_p(r):0\leq r \leq
[\frac{m}{p^{\nu_p(n)}(p-1)}] \}.$$

If $p=2,$ then
$$\theta_p(y_m)=(1+\frac{1}{4}y_m)^{\frac{1}{2}}.$$
So,
$$\psi^n(\theta_p(y_m))=(1+\frac{1}{4}\psi^n(y_m))^{\frac{1}{2}}.$$
Let $n=p_{1}^{r_1} \ldots p_{s}^{r_s} 2^d .$ Then
$1+\frac{1}{4}\psi^n(y_m)=\sum_{j=1}^{n}c_{j}y_{m}^{j}$
with $\nu_p(c_j)=-2$ if $j$ is a multiple of
$2^d,$ $\nu_p(c_j)\geq 0$ for $j<2^d,$ and
$\nu_p(c_j)\geq -1$ for all other values of $j.$
Hence, $1+\frac{1}{4}\psi^n(y_m)$ is a strict
D-series at $2$ of type $(\alpha,\beta)$ where $$
\alpha_{p^{\prime}}=\left\{ \begin{array}{ll} 2^d
& p^{\prime}=2 \\ \infty & p^{\prime} \neq 2
\end{array} \right.  \hspace{2cm}, \,
\beta_{p^{\prime}}=\left\{ \begin{array} {ll} 2 &
p^{\prime}=2 \\ \infty & p^{\prime} \neq 2.
\end{array}\right.  $$ So,
$$\nu_2(e_t(1+\frac{1}{4}\psi^n(y_m)))=\max\{2r+\nu_2(r):0\leq
r \leq [\frac{t}{2^d}] \}.$$ Hence,
$$\nu_2(b_m(r\xi_m({\Bbb
C})^n-2))=\nu_2(e_t((1+\frac{1}{4}
\psi^n(y_m)))^{\frac{1}{2}})$$
$$=\max\{2r+\nu_2(2r):0\leq r \leq [\frac{t}{2^d}]
\}=\max\{r+\nu_2(r):0\leq r \leq [\frac{m}{2^d}]
\}.$$ So, we have:
\begin{thm} Let $n \in {\Bbb N}$ and $p$ be any
prime number then $$ \nu_p(b_m(r\xi_m({\Bbb
C})^n-2))=\max\{r+\nu_p(r):0\leq r \leq
[\frac{m}{p^{\nu_p(n)}(p-1)}] \}.$$ \end{thm}
REMARK.  A similar proof of Theorem~{3.2} when $n$
is a power of a prime $p$ has been obtained
independently by Diba\v{g} \cite{dib2}.

$\begin{array}{c}\mbox{University of California }
\\ \mbox{ Department of Mathematics }\\ \mbox{
Berkeley, CA 94720, USA } \\ \mbox{ e-mail :
obiedat@math.berkeley.edu } \end{array}$
\end{document}